\documentclass[12pt]{article}
\usepackage{amsmath}
\usepackage{amsfonts}
\usepackage{mathrsfs}
\usepackage{amssymb}
\usepackage[usenames]{color}
\usepackage[numbers,sort&compress]{natbib}
\usepackage{subfigure}

\def\cl{\centerline}
\def\ni{\noindent}
\def\vs{\vspace*}

\def\a{\alpha}
\def\b{\beta}
\def\D{\Delta}

\def\la{\lambda}
\def\pa{\partial}

\def\CVir{{\frak {Vir}}}
\def\CB{{\frak {B}}(p)}
\def\CBB{{\frak {B}}}
\def\BB{{\cal B}(q)}
\def\sp{{\rm{span}}}

\def\Z{\mathbb{Z}}

\def\C{\mathbb{C}}
\def\Q{\mathbb{Q}}
\def\QED{\hfill$\Box$}

\numberwithin{equation}{section}
\newtheorem{theo}{Theorem}[section]
\newtheorem{defi}[theo]{Definition}
\newtheorem{coro}[theo]{Corollary}
\newtheorem{lemm}[theo]{Lemma}
\newtheorem{rema}[theo]{Remark}
\newtheorem{prop}[theo]{Proposition}
\newtheorem{clai}{Claim}

\begin{document}
\begin{center}
{\bf\Large Classification of finite irreducible conformal modules over a class of Lie conformal algebras of Block type}
\footnote{
$^{\dag}$Corresponding author: chgxia@cumt.edu.cn (C.~Xia).}
\end{center}

\cl{Yucai Su$^{\,*}$, \ Chunguang Xia$^{\,\dag}$, \ Lamei Yuan$^{\,\ddag}$}

\cl{\small $^{\,*}$School of Mathematical Sciences, Tongji University, Shanghai 200092, China}

\cl{\small $^{\,\dag}$School of Mathematics, China University of Mining and Technology, Xuzhou 221116, China}

\cl{\small $^{\,\ddag}$Academy of Fundamental and Interdisciplinary Sciences,\vspace{-2pt}}
\cl{\small Harbin Institute of Technology, Harbin 150080, China}

\cl{\small Email: ycsu@tongji.edu.cn, chgxia@cumt.edu.cn, lmyuan@hit.edu.cn}
\vspace{5pt}

\footnotesize
\noindent{{\bf Abstract.}
We classify finite irreducible conformal modules over
a class of infinite Lie conformal algebras $\CB$ of Block type, where $p$ is a nonzero complex number.
In particular, we obtain that a finite irreducible conformal module over $\CB$ may be a nontrivial
extension of a finite conformal module over $\CVir$ if $p=-1$, where $\CVir$ is a Virasoro conformal subalgebra of $\CB$.
As a byproduct, we also obtain the classification of finite irreducible conformal modules over a series of finite Lie conformal algebras ${\frak b}(n)$ for $n\ge1$.
\vs{5pt}

\ni{\bf Key words:} finite conformal module; Lie conformal algebras of Block type; Virasoro conformal algebra

\ni{\it Mathematics Subject Classification (2010):} 17B10; 17B65; 17B68; 17B69.}

\small
\section{Introduction}
Lie conformal algebras, introduced by Kac \cite{K1}, encode the singular part of the operator
product expansion of chiral fields in conformal field theory.
%It is in some sense a generalization of a Lie algebra in that it is also a ``Lie algebra'' but in a different pseudo-tensor category.
%It turned out to be an adequate tool for the study of infinite dimensional Lie algebras satisfying the locality property.
%Conformal modules are basic tools for the construction of free field realization of infinite dimensional Lie (super)algebras in conformal field theory.
%In recent years, the structure and representation theory of Lie conformal algebras have been extensively studied.
The theory of finite Lie conformal algebras has been greatly developed in the last two decades (e.g., \cite{BKV,CK,CKW,DK,DeK2,K1,K2,Z1}).
Finite simple Lie conformal algebras were classified in \cite{DK},  which shows that a finite simple Lie conformal algebra is isomorphic to either
the Virasoro conformal algebra or a current conformal algebra $\rm Cur\,\mathfrak{g}$ over a
 simple finite-dimensional Lie algebra $\mathfrak{g}$.
The theory of conformal modules and their extensions was developed in \cite{CK,CKW},
and the cohomology theory was developed in \cite{BKV} and further  in \cite{DeK2}.
For super cases, the structure and representation theories have also been developed in recent years, see \cite{FK,FKR,BKL3,MZ} and the references therein.

However, the theory of infinite Lie conformal algebras is far from being well developed.
The most important example %of infinite simple Lie conformal algebra
is the general Lie conformal algebra $gc_N$,
%any module $M=\C[\partial]^N$ over a Lie conformal algebra $R$ is obtained via a homomorphism $R\rightarrow gc_N$.
%In other words, the general Lie conformal algebra $gc_N$
which plays the same role in the theory of
Lie conformal algebras as the general Lie algebra $gl_N$ does in the theory of Lie algebras.
The important difference between usual and conformal algebras is that $gc_N$ is infinite.
Due to these reasons, the general Lie conformal algebra $gc_N$ and its subalgebras have been studied by many authors (e.g., \cite{BKL1,BKL2,BM,DK,DeK1,S1,SYue,Z2}).
Recently, some interesting examples of infinite Lie conformal algebras associated with infinite-dimensional loop Lie algebras were constructed and studied (e.g., \cite{FSW,FSX,WCY}).

In this paper, we focus on another class of infinite Lie conformal algebras $\CB$
with $p$ being a nonzero complex number, where $\CB$ has a $\C[\partial]$-basis $\{L_i\,|\,i\in\Z_+\}$
and $\lambda$-brackets
\begin{equation}\label{brackets}
[L_i\, {}_\lambda \, L_j] = ((i+p)\partial+(i+j+2p)\lambda) L_{i+j}.
\end{equation}
We refer to $\CB$'s as {\it Lie conformal algebras of Block type} due to their
relations with some Lie algebras of Block type (see Remark~\ref{remark-name-reason}).
There are some interesting features on this class of Lie conformal algebras.
Firstly, each $\CB$ contains a Virasoro conformal subalgebra.
Set $L=\frac{1}{p}L_0\in\CB$. By \eqref{brackets}, we see that
$[L\, {}_\lambda \, L]=(\partial+2\lambda) L.$ Namely, the subalgebra
\begin{equation}\label{equ-CVir}
\CVir=\C[\partial]L %=\C[\partial]L_0
\end{equation}
of $\CB$ is exactly the Virasoro conformal algebra.
Secondly, the special case $\CBB(1)$ has close relation with the general Lie conformal algebra $gc_1$. In fact, $\CBB(1)$ is a maximal subalgebra of
the associated graded conformal algebra ${\rm{gr}}\,gc_1$ of the filtered algebra $gc_1$ \cite{SYue}.
Thirdly, there are embedding relations among $\CB$'s. For any integer $n\ge 1$, $\CB$ can be embedded into $\CBB(np)$ via $L_i\mapsto \frac{1}{n}L'_{ni}$.
%In particular, the special case $\CBB(1)$ can be embedded into $\CBB(n)$ via $L_i\mapsto \frac{1}{n}L'_{ni}$.
Finally, $\CBB(-n)$ contains a series of finite Lie conformal quotient algebras (cf.~\eqref{equ-quotiont-special})
$$
{\frak b}(n)=\CBB(-n)/\CBB(-n)_{\langle n+1 \rangle},% \ \ \text{with}\ \ n\ge 1,
$$
including the Heisenberg-Virasoro conformal algebra ${\frak b}(1)$
and Schr$\ddot{\rm o}$dinger-Virasoro conformal algebra ${\frak b}(2)$ as the first two cases (see Subsection~2.2).
Due to these observations, it seems to be interesting for us to study the Lie conformal algebra $\CB$ systematically.

%There are some interesting relations on the Lie conformal algebra $\CB$:
%\begin{itemize}
%  \item[(R1)] {\bf(Relation with the Virasoro conformal algebra)}~~Set $L=\frac{1}{p}L_0\in\CB$. By \eqref{brackets}, we see that
%  $$[L\, {}_\lambda \, L]=(\partial+2\lambda) L.$$ Namely, the subalgebra $\CVir=\C[\partial]L_0$ of $\CB$ is
%  exactly the Virasoro conformal algebra.
%
%  \item[(R2)] {\bf(Relation with the general Lie conformal algebras)}~~The special case $\CBB(1)$ is a maximal subalgebra of
%  the associated graded conformal algebra ${\rm{gr}}\,gc_1$ of the filtered algebra $gc_1$ \cite{SYue}.
%
%  \item[(R3)] {\bf(Relation among themselves)}~~For any integer $n\ge 1$, $\CB$ can be embedded into $\CBB(np)$ via $L_i\mapsto \frac{1}{n}L'_{ni}$.
%  In particular, the special case $\CBB(1)$ can be embedded into $\CBB(n)$ via $L_i\mapsto \frac{1}{n}L'_{ni}$.
%\end{itemize}

One of the most important problems regarding infinite Lie conformal algebras is the classification of
finite irreducible conformal modules (FICMs). This problem for the general Lie conformal algebra $gc_N$ was solved by  Kac, Radul and Wakimoto, see also \cite{BM,K2}.
In this paper, we consider this problem for $\CB$.
Clearly, any conformal module over $\CVir\subset\CB$ can be trivially extended to a conformal module over the whole $\CB$.
Our main results indicate that a FICM over $\CB$ may be a nontrivial
extension of a finite conformal module over $\CVir$ if $p=-1$ (see Table~1).
As a byproduct, we also obtain the classification of FICMs over the finite Lie conformal algebra ${\frak b}(n)$ (see Table~2).

\begin{table}[h]
\centering\small
\subtable[Nontrivial FICMs over $\CB$]{
\begin{tabular}{c|c|c}
\hline
 $\CB$  & FICMs & Reference \\
\hline
$p\ne -1$  & $M_{\D,\a}$ & Theorems~\ref{thm-rank-one-modules} and \ref{thm-classification} \\
\hline
$p=-1$ & $M_{\D,\a,\b}$ & Theorems~\ref{thm-rank-one-modules} and \ref{thm-classification}\\
\hline
\end{tabular}}
\quad\quad
\subtable[Nontrivial FICMs over ${\frak b}(n)$]{
\begin{tabular}{c|c|c}
\hline
${\frak b}(n)$  & FICMs & Reference \\
\hline
$n> 1$  & $M_{\D,\a}$ & Corollary~\ref{classification-b} \\
\hline
$n=1$ & $M_{\D,\a,\b}$ & Corollary~\ref{classification-b}\\
\hline
\end{tabular}}
\end{table}

%\begin{center}
%~~\mbox{Table 1: Nontrivial FICMs over $\CB$} \quad\quad\quad~   \mbox{Table 2: Nontrivial FICMs over ${\frak b}(n)$}
%\end{center}

Our main results imply that $\CB$'s are examples of infinite conformal algebras that have
a finite module. The classification of such infinite conformal algebras (even for simple ones) is a challenging problem as Kac stated in \cite{K3}.
A more general problem is to classify all simple conformal algebras of finite growth.
In fact, the ``bigger'' $\CB$'s are such examples, which can be obtained by replacing the index set $\Z_+$ of
the $\C[\partial]$-base elements $L_i$ of $\CB$'s by $\Z$ (see Remark~\ref{remark-bigger-alg}).
% Our proof uses some ideas of \cite{SYue} and \cite{CK}.

We also would like to point out that our techniques used here may be  applied to
analogous problems of (both infinite and finite) Lie conformal algebras which are closely related to the Lie conformal algebra $\CB$. This is also our
motivation for writing this paper.

This paper is arranged as follows.
In Section~2, we recall some definitions on Lie conformal algebras and present some
related algebraic structures of $\CB$. %the Lie conformal algebra $\CB$.
Then, in Section~3, we give three technical lemmas,
which will gradually reduce our classification problem.
In Section~4, we first consider the problem for the rank one case without the irreducibility assumption.
%(see Theorem~\ref{rank-one-modules}).
Finally, in Section~5, we completely classify FICMs over $\CB$ by showing that they must be free of rank one.
By the last observation mentioned above, we also obtain the classification of FICMs over the finite Lie conformal algebra ${\frak b}(n)$.

\section{Preliminaries}

\subsection{Definitions}
First we list some definitions that will be used. All of them are collected or reorganized from \cite{DK,K1,K2,K3}.

\begin{defi}\label{def-conforma-algebra}\rm
A {\it Lie conformal algebra} $R$ is a $\C[\partial]$-module endowed with a
$\C$-linear map $R\otimes R\rightarrow \C[\lambda]\otimes R$, $a\otimes b\rightarrow [a\,{}_\lambda \,b]$
called $\lambda$-bracket, and satisfying the following axioms ($a,\,b,\,c\in R$):
\begin{equation*}
\aligned
\mbox{(conformal sesquilinearity)}&~~~~[\partial a\,{}_\lambda \,b]=-\lambda[a\,{}_\lambda\, b],\ \ \ \
[a\,{}_\lambda \,\partial b]=(\partial+\lambda)[a\,{}_\lambda\, b],\\
\mbox{(skew-symmetry)}&~~~~[a\, {}_\lambda\, b]=-[b\,{}_{-\lambda-\partial}\,a],\\
\mbox{(Jacobi identity)}&~~~~[a\,{}_\lambda\,[b\,{}_\mu\, c]]=[[a\,{}_\lambda\, b]\,{}_{\lambda+\mu}\, c]+[b\,{}_\mu\,[a\,{}_\lambda \,c]].
\endaligned
\end{equation*}
\end{defi}

\begin{defi}\label{def-conformal-module}\rm
A {\it conformal module} $M$ over a Lie conformal algebra $R$ is a $\C[\partial]$-module endowed with a $\lambda$-action
$R\otimes M\rightarrow \C[\lambda]\otimes M$, $a\otimes v\rightarrow a\,{}_\lambda \,b$, such that ($a,\,b\in R$, $v\in M$)
\begin{equation*}
(\partial a)\,{}_\lambda\, v=-\lambda a\,{}_\lambda\, v,\ \ \ \ \ a{}\,{}_\lambda\, (\partial v)=(\partial+\lambda)a\,{}_\lambda\, v,\ \ \ \
[a\,{}_\lambda\, b]\,{}_{\lambda+\mu}\, v = a\,{}_\lambda\, (b{}\,_\mu\, v)-b\,{}_\mu\,(a\,{}_\lambda\, v).
\end{equation*}
A conformal $R$-module $M$ is called {\it finite} if it is finitely generated over $\C[\partial]$.
\end{defi}

\begin{defi}\label{def-annihilation-algebra}\rm
An {\it  annihilation algebra} ${\cal A}(R)$ of a Lie conformal algebra $R$ is a Lie algebra with
$\C$-basis $\{a_{(n)}\,|\,a\in R,\,n\in\Z_+\}$ and relations
$$
[a_{(m)}, b_{(n)}]=\sum_{k\in\Z_+}{\binom m k}(a_{(k)}b)_{(m+n-k)}, \ \ \ \ \ (\partial a)_{(n)}=-n a_{(n-1)},
$$
where $a_{(k)}b$ is called the $k$-product, given by
$[a\,{}_\lambda \, b]=\sum_{k\in\Z_+}\lambda^{(k)}a_{(k)}b$ with $\lambda^{(k)}=\frac{\lambda^k}{k!}$.
Furthermore, an {\it extended annihilation algebra} ${\cal A}(R)^e$ of $R$ is defined by
${\cal A}(R)^e=\C T\ltimes{\cal A}(R)$ with $[T, a_{(n)}]=-n a_{(n-1)}$.
\end{defi}

\begin{prop}\label{prop-observation}
A conformal module $M$ over a Lie conformal algebra $R$ is the same as
a module over the Lie algebra ${\cal A}(R)^e$ satisfying $a_{(n)}v=0$ for $a\in R$, $v\in M$, $n\gg 0$.
\end{prop}

\subsection{Quotient algebras of $\CB$}

Start from $\CB$, one can obtain many interesting finite Lie conformal algebras.
Consider, for example, the quotient algebras $\CB_{[n]}$ with $n\in\Z_+$, defined by
\begin{equation}\label{equ-quotiont}
\CB_{[n]}=\CB/\CB_{\langle n+1\rangle}, \quad \mbox{where} \quad \CB_{\langle n\rangle}= {\oplus}_{i\ge n}\C[\partial]L_i.
\end{equation}
Note that $\CB_{[0]}$ is isomorphic to the Virasoro conformal algebra $\CVir$. All other $\CB_{[n]}$ with $n\ge 1$ are non-simple.
The special quotients
\begin{equation}\label{equ-quotiont-special}
{\frak b}(n)=\CBB(-n)_{[n]}=\CBB(-n)/\CBB(-n)_{\langle n+1 \rangle} \ \ \text{with}\ \ n\ge 1,
\end{equation}
seem to be particularly interesting.  Let us look at the first two cases.
%(here, we still use $L_i$ to denote the $\C[\partial]$-base elements of ${\frak b}(n)$):
\begin{itemize}
  \item[(Q1)] Case $n=1$. Set $L=-\bar{L}_0$ and $M=\bar{L}_1$. By \eqref{brackets} and \eqref{equ-quotiont-special}, we see that
  $$[L\, {}_\lambda \, L]=(\partial+2\lambda) L,\ \ [L\, {}_\lambda \, M]=(\partial+\lambda) L, \ \
  [M\, {}_\lambda \, L]=\lambda L,\  \ [M\, {}_\lambda \, M]=0.$$
  Namely, ${\frak b}(1)=\C[\partial]L\oplus\C[\partial]M$ is exactly the Heisenberg-Virasoro conformal algebra \cite{SY},
  which is associated with the Heisenberg-Virasoro Lie algebra \cite{ADKP}.

  \item[(Q2)] Case $n=2$. Set $L=-\frac{1}{2}\bar{L}_0$, $Y=\bar{L}_1$ and $M=-\bar{L}_2$. By \eqref{brackets} and \eqref{equ-quotiont-special},
  we see that  (other components vanish)
$$
\begin{aligned}
&  [L\, {}_\lambda \, L]=(\partial+2\lambda) L, && [L\, {}_\lambda \, Y]=(\partial+\frac{3}{2}\lambda) Y, && [L\, {}_\lambda \, M]=(\partial+\lambda) M,\\[-2pt]
&  [Y\, {}_\lambda \, L]=(\frac{1}{2}\partial+\frac{3}{2}\lambda) Y, && [Y\, {}_\lambda \, Y]=(\partial+2\lambda) M, && [M\, {}_\lambda \, L]=\lambda M.
\end{aligned}
$$
  Namely, ${\frak b}(2)=\C[\partial]L\oplus\C[\partial]Y\oplus\C[\partial]M$ is exactly the Schr$\ddot{\rm o}$dinger-Virasoro conformal algebra \cite{SY},
  which is associated with the Schr$\ddot{\rm o}$dinger-Virasoro Lie algebra \cite{He}.
\end{itemize}

\subsection{Annihilation algebra of $\CB$}

Denote by ${\cal A}(\CB)$ and ${\cal A}(\CB)^e$ the annihilation algebra and extended annihilation algebra of $\CB$, respectively.
Their concrete Lie structures are as follows.

\begin{lemm}\label{annihilation-algebra-str}
\baselineskip1pt\lineskip7pt\parskip-1pt
\begin{itemize}\parskip-1pt
  \item[{\rm(1)}] The annihilation algebra ${\cal A}(\CB)$ is given by
  $$
  {\cal A}(\CB)=\sp_{\C}\{L_{i,m}\,|\,i\in\Z_{+},\,m\in\Z_{\ge-1}\}
  $$ with relations
$$[L_{i,m},L_{j,n}]=((j+p)(m+1)-(i+p)(n+1)) L_{i+j,m+n}.$$
  \item[{\rm(2)}] The extended annihilation algebra ${\cal A}(\CB)^e$ is given by
  $$
  {\cal A}(\CB)^e=\sp_{\C}\{L_{i,m},\,T\,|\,i\in\Z_{+},\,m\in\Z_{\ge-1}\}
  $$ with relations as in {\rm (1)} and $[T, L_{i,m}]=-(m+1)L_{i,m-1}.$
\end{itemize}
\end{lemm}

\ni{\it Proof.}\ \
For $i,j\in\Z_+$, by \eqref{brackets} and Definition~\ref{def-annihilation-algebra}, we have
$$
((i+p)\partial+(i+j+2p)\lambda) L_{i+j}=[L_i\, {}_\lambda \, L_j]=\sum_{k\in\Z_+}\lambda^{(k)}{L_i}_{(k)}{L_j}
={L_i}_{(0)}{L_j}+\lambda{L_i}_{(1)}{L_j}+\sum_{k\ge 2}\lambda^{(k)}{L_i}_{(k)}{L_j},
$$
which implies that
$$
{L_i}_{(0)}{L_j}=(i+p)\partial L_{i+j}, \ \ {L_i}_{(1)}{L_j}=(i+j+2p)L_{i+j}, \ \ {L_i}_{(k)}{L_j}=0\ \ \mbox{for}\ \ k\ge2.
$$
Then, for $m,n\in\Z_+$, by Definition~\ref{def-annihilation-algebra}, we have
\begin{eqnarray}
% \nonumber to remove numbering (before each equation)
\nonumber[(L_i)_{(m)}, (L_j)_{(n)}] &\!\!=\!\!& \sum_{k\in\Z_+}{\binom m k}({L_i}_{(k)}{L_j})_{(m+n-k)}\\
\nonumber  &\!\!\!=\!\!\!& ({L_i}_{(0)}{L_j})_{(m+n)}+m({L_i}_{(1)}{L_j})_{(m+n-1)} \\
\nonumber  &\!\!\!=\!\!\!& -(m+n)(i+p)(L_{i+j})_{(m+n-1)}+m(i+j+2p)(L_{i+j})_{(m+n-1)} \\
\nonumber  &\!\!\!=\!\!\!& (m(j+p)-n(i+p))(L_{i+j})_{(m+n-1)},
\end{eqnarray}
and $[T,(L_i)_{(m)}]=-m (L_i)_{(m-1)}$. Hence, this lemma holds by setting $L_{i,m}=(L_i)_{(m+1)}$ for $m\ge-1$.
\QED

\begin{rema}\label{remark-quotionts}\rm
Denote by ${\cal A}(\CB_{[n]})$ and ${\cal A}(\CB_{[n]})^e$ the annihilation algebra and extended annihilation algebra of the quotient algebra $\CB_{[n]}$ of $\CB$, respectively.
Then %they respectively have $\C$-bases %(cf.~\eqref{equ-quotiont})
\begin{eqnarray}
% \nonumber to remove numbering (before each equation)
\nonumber {\cal A}(\CB_{[n]}) &\!\!=\!\!& \sp_{\C}\{\bar{L}_{i,m}\,|\,0\le i\le n,\,m\in\Z_{\ge-1}\},\\
\nonumber {\cal A}(\CB_{[n]})^e &\!\!\!=\!\!\!& \sp_{\C}\{\bar{L}_{i,m},\,T\,|\,0\le i\le n,\,m\in\Z_{\ge-1}\},
\end{eqnarray}
%$$\{\bar{L}_{i,m}\,|\,0\le i\le n,\,m\in\Z_{\ge-1}\}\ \ \mbox{and}\ \ \{\bar{L}_{i,m},\,T\,|\,0\le i\le n,\,m\in\Z_{\ge-1}\},$$
which satisfy the same relations as in Lemma~\ref{annihilation-algebra-str}
(of course, $\bar{L}_{i,m}$ with $i>n$ will be viewed as zero if it technically appears).
Clearly, ${\cal A}(\CB_{[n]})$ and ${\cal A}(\CB_{[n]})^e$
can be also viewed as quotient algebras of ${\cal A}(\CB)$ and ${\cal A}(\CB)^e$, respectively.
%Here, all elements $\bar{L}_{i,m}$ with $i> n$ are conventionally viewed as zero if they technically appear in the context.
\end{rema}

\begin{rema}\label{remark-name-reason}\rm
The annihilation algebra ${\cal A}(\CB)$ has close relation with the Lie algebra $\BB$ of Block type studied in \cite{SXX1,SXX2,XZ}.
Consider a two-parameter Lie algebra ${\cal B}(p,q)$ of Block type with $p,q\in\C$, where ${\cal B}(p,q)$ has a $\C$-basis
$\{L_{i,m}\,|\,i,\,m\in\Z\}$ and relations
$$[L_{i,m},L_{j,n}]=((j+p)(m+q)-(i+p)(n+q)) L_{i+j,m+n}.$$
Then, the special case ${\cal B}(p,q)$ with $q=1$ contains ${\cal A}(\CB)$ as a subalgebra, while the special case ${\cal B}(p,q)$ with $p=0$ contains $\BB$ as a subalgebra.
Hence, we refer to $\CB$'s as {\it Lie conformal algebras of Block type}.
\end{rema}

\section{Three technical lemmas} %reduction lemmas

In this section, we shall give three lemmas for later use.

Our first lemma is the following result, which will reduce the classification of finite conformal modules over $\CB$
to the classification of such modules over a finite conformal quotient algebra of $\CB$.
%It is well-known that a module over a quotient algebra of an arbitrary algebra $\frak{A}$ can be canonically viewed as a module over $\frak{A}$,
%but the converse is usually not true. The situation seems different for $\CB$.
%We have the following result.
The case $p=1$ is in fact implicit in \cite{SYue}.
Our proof here is a straightforward generalization of that in \cite{SYue}, but with a slightly different discussion (at the end of Case~1).
%Here, we still include the proof for completeness.

\begin{lemm}\label{key-lemma-1}
Let $M$ be a nontrivial finite conformal module over $\CB$. Then the $\lambda$-action of $L_i\in\CB$ on $M$ is trivial for all $i\gg 0$.
In particular, a finite conformal module over $\CB$ is simply a finite conformal module over $\CB_{[n]}$ for some big enough integer $n$,
where $\CB_{[n]}$ is defined by \eqref{equ-quotiont}.
\end{lemm}

\ni{\it Proof.}\ \
First, recall that \cite{CK} a finite irreducible conformal module over
$\CVir=\C[\partial]L$ is isomorphic to either a free conformal module of rank one $M_{\D,\a}=\C[\partial]v$ with an action defined by
\begin{equation}\label{CVir-actions}
L\,{}_\lambda\, v=(\partial+\D\lambda+\a)v
\end{equation}
for some $\D,\a\in\C$ with $\D\ne0$, or a one-dimensional trivial module $\C c_{\a}$ with an action defined by
%\begin{equation}\label{trivial-action}
%L\,{}_\lambda\, c_\a=0, \quad \partial c_{\a}=\a c_\a
%\end{equation}
$L\,{}_\lambda\, c_\a=0$, $\partial c_{\a}=\a c_\a$
for some $\a\in\C$.

By regarding $M$ as a module over $\CVir\subset\CB$, we can choose a composition series
$$
M=M_N\supset M_{N-1}\supset\cdots\supset M_1\supset M_0=0,
$$
such that for each $1\le k\le N$, the composition factor $\overline{M}_k=M_k/M_{k-1}$ is either a free conformal module of rank one
$M_{\D_k,\a_k}$ with $\D_k\ne 0$, or a one-dimensional trivial module $\C c_{\a_k}$.
Denote by $\bar{v}_k$ a $\C[\partial]$-generator of $\overline{M}_k$ and $v_k\in M_k$ the preimage of $\bar{v}_k$.
Then $\{v_k\,|\,1\le k\le N\}$ is a $\C[\partial]$-generating set of $M$, and the $\lambda$-action of $L_0$ on
$v_k$ is a $\C[\partial,\lambda]$-combination of $v_1, \ldots, v_k$.

We claim that the $\lambda$-action of $L_i$ on $v_1$ is trivial for all $i\gg 0$. Namely,
\begin{equation}\label{claim-trivial}
L_i\,{}_\lambda\, v_1=0\ \ \ \mbox{for all}\ \ \ i\gg 0.
\end{equation}
%\begin{clai}
%The $\lambda$-action of $L_i$ on $v_1$ is trivial for all $i\gg 0$. Namely, $L_i\,{}_\lambda\, v_1=0$ for all $i\gg 0$.
%\end{clai}
Fix $i\gg 0$ and assume that $L_i\,{}_\lambda\, v_1\ne 0$.
Let $k_i\ge 1$ be the largest integer such that $L_i\, {}_\lambda\, v_1 \notin M_{k_i-1}$.
We proceed to derive a contradiction. We only need to consider the following four cases.

\vskip2pt
{\bf Case 1:} $M_1=M_{\D_1,\a_1}$, $\overline{M}_{k_i}=M_{\D_{k_i},\a_{k_i}}$.
\vskip2pt

By assumption, we can write
\begin{equation}\label{case1}
L_i\, {}_\lambda\, v_1 \equiv f_i(\partial,\lambda)v_{k_i} (\mbox{mod}\ M_{k_i-1})\ \ \mbox{for some}\ \ 0\ne f_i(\partial,\lambda)\in\C[\partial,\lambda].
\end{equation}
Considering the action of  the operator $L_0\,{}_\mu$ on \eqref{case1},  by Definition~\ref{def-conformal-module}, we obtain (note that $L_0=pL$)
\begin{equation}\label{case1-L0-action}
p(\partial+\D_{k_i}\mu+\a_{k_i})f_i(\partial+\mu,\lambda)=((i+p)\mu-p\lambda)f_i(\pa,\mu+\la)+p(\pa+\lambda+\D_1\mu+\a_1)f_i(\pa,\lambda).
\end{equation}
In particular, taking $\pa=0$ in \eqref{case1-L0-action}, we have %(note that $\D_{k_i}\ne0$)
\begin{equation}\label{case1-L0-action-0}
f_i(\mu,\lambda)=\frac{1}{p(\D_{k_i}\mu+\a_{k_i})}(((i+p)\mu-p\lambda)f_i(0,\mu+\lambda)+p(\lambda+\D_1\mu+\a_1)f_i(0,\lambda)).
\end{equation}
Taking $\pa=-\D_{k_i}\mu-\a_{k_i}$ and $\lambda=(1+\frac{i}{p})\mu$ in \eqref{case1-L0-action}, and then using \eqref{case1-L0-action-0}, we obtain
$$
(1+\frac{i}{p})((1+\D_{k_i})\mu+\a_{k_i})f_i(0,(1+\frac{i}{p}-\D_{k_i})\mu-\a_{k_i})
=((1+\frac{i}{p}-\D_1\D_{k_i})\mu-\D_1\a_{k_i}+\a_1)f_i(0,(1+\frac{i}{p})\mu).
$$
Denote by $m_i$ the degree of $f_i(0,\lambda)$.  Equating the coefficients of $\mu^{m_i+1}$ in the above equation, we obtain
\begin{equation}\label{case1-final}
(1+\D_{k_i})(1+\frac{i}{p})(1+\frac{i}{p}-\D_{k_i})^{m_i}
=(1+\frac{i}{p}-\D_1\D_{k_i})(1+\frac{i}{p})^{m_i}.
\end{equation}
Set $q(i)=1+\frac{i}{p}$. Then either $q(i)\gg 0$ or $q(i)\ll 0$ when $i\gg0$.
Comparing the coefficients of $q(i)^{m_i+1}$ in \eqref{case1-final},
one can see that it cannot hold for any $m_i\ge 0$ (note that $\D_{k_i}\ne0$), a contradiction.

\vskip2pt
{\bf Case 2:} $M_1=\C c_{\a_1}$, $\overline{M}_{k_i}=M_{\D_{k_i},\a_{k_i}}$.
\vskip2pt

As in Case~1, we can still assume \eqref{case1}. Considering the action of  the operator $L_0\,{}_\mu$ on \eqref{case1},
by Definition~\ref{def-conformal-module}, we obtain
\begin{equation}\label{case2-L0-action}
p(\partial+\D_{k_i}\mu+\a_{k_i})f_i(\partial+\mu,\lambda)=((i+p)\mu-p\lambda)f_i(\pa,\mu+\lambda).
\end{equation}
Taking $\pa=\mu=0$ in \eqref{case2-L0-action}, we obtain $f_i(0,\lambda)=0$.
Using this in \eqref{case2-L0-action} with $\pa=0$, we obtain $f_i(\mu,\lambda)=0$, a contradiction.

\vskip2pt
{\bf Case 3:} $M_1=M_{\D_{1},\a_{1}}$, $\overline{M}_{k_i}=\C c_{\a_{k_i}}$.
\vskip2pt

In this case, since $\pa$ acts on $\bar{v}_{k_i}$ as the scalar $\a_{k_i}$, we can write
\begin{equation}\label{case3}
L_i\, {}_\lambda\, v_1 \equiv f_i(\lambda)v_{k_i} (\mbox{mod}\ M_{k_i-1})\ \ \mbox{for some}\ \ 0\ne f_i(\lambda)\in\C[\lambda].
\end{equation}
Considering the action of  the operator $L_0\,{}_\mu$ on \eqref{case3},  by Definition~\ref{def-conformal-module},
we obtain
\begin{equation}\label{case3-L0-action}
0=((i+p)\mu-p\lambda)f_i(\mu+\lambda)+p(\pa+\lambda+\D_1\mu+\a_1)f_i(\lambda).
\end{equation}
Equating the coefficients of $\pa$ in \eqref{case3-L0-action}, we immediately obtain $f_i(\lambda)=0$, a contradiction.

\vskip2pt
{\bf Case 4:} $M_1=\C c_{\a_1}$, $\overline{M}_{k_i}=\C c_{\a_{k_i}}$.
\vskip2pt

As in Case~3, one can easily obtain $f_i(\lambda)=0$, a contradiction.

Now, start from \eqref{claim-trivial}, one can inductively show that $L_i\,{}_\lambda\, v_k=0$ for $1\le k\le N$.
Hence, the $\lambda$-action of $L_i$ on $M$ is trivial. This completes the proof.
\QED

\begin{rema}\label{remark-bigger-alg}\rm
Clearly, $\CB$ is non-simple. If we replace the index set $\Z_+$ of
the $\C[\partial]$-base elements $L_i\in\CB$ by $\Z$, we obtain a simple Lie conformal algebra.
By Lemma~\ref{key-lemma-1} and the simplicity of this ``bigger'' $\CB$, one can easily show that
it has no nontrivial finite conformal module.
In particular, as pointed out in \cite{SYue}, this provides a class of examples of finitely freely generated simple Lie conformal algebras of linear
growth that cannot be embedded into $gc_N$ for any $N$.
\end{rema}

Our second lemma is a result concerning the representations of certain subquotient algebra of the annihilation algebra ${\cal A}(\CB)$.
Let $k\ge 0$ and $N\ge 0$ be two fixed integers. Define a Lie algebra ${\cal{G}}_{k,N}$ by
$$
{\cal{G}}_{k,N}=\sp_{\C}\{J_{i,m}\,|\,0\le i\le k,\,0\le m\le N\}
$$
with relations
\begin{equation*}
[J_{i,m},J_{j,n}]=
\left\{\begin{array}{ll}
((j+p)(m+1)-(i+p)(n+1)) J_{i+j,\,m+n}&\mbox{if \ }i+j\le k,\,m+n\le N,\\[6pt]
0&\mbox{otherwise.}\end{array}\right.
\end{equation*}
%$[J_{i,m},J_{j,n}]=((j+p)(m+1)-(i+p)(n+1)) J_{i+j,\,m+n}$.
%Here, $J_{i,m}$ with $i>k$ or $m>N$ will be viewed as zero if it technically appears.
By Lemma~\ref{annihilation-algebra-str}, ${\cal{G}}_{k,N}$ can be viewed as a subquotient algebra of ${\cal A}(\CB)$.

\begin{lemm}\label{key-lemma-3}
Let $V$ be a nontrivial finite dimensional irreducible module over ${\cal{G}}_{k,N}$. Then $\dim V=1$.
\end{lemm}

\ni{\it Proof.}\ \
Denote by $\Q_+$ the set of all positive rational numbers. We divide the proof into the following two cases.

\vskip2pt
{\bf Case 1:} $p\notin\Q_+$.
\vskip2pt

Consider the following decomposition of ${\cal{G}}_{k,N}$:
$$
{\cal{G}}_{k,N}=\C J_{0,0}+{\cal{K}}, \mbox{\ \ where\ \ }{\cal{K}}={\cal{G}}_{k,N}\backslash\C J_{0,0}.
$$
Clearly, ${\cal{K}}$ is a nilpotent ideal of ${\cal{G}}_{k,N}$. For any $J_{i,m}\in{\cal{K}}$,
we have $[J_{0,0},J_{i,m}]=(i-pm) J_{i,\,m}$, and $i-pm\ne0$ since $p\notin\Q_+$.
It follows that ${\cal{K}}$ is a completely reducible $\C J_{0,0}$-module with no trivial summand.
By \cite[Lemma~1]{CK}, ${\cal{K}}$ acts trivially on $V$. Hence, $V$ can be viewed as a finite dimensional module
over $\C J_{0,0}$, and so $\dim V=1$.

\vskip2pt
{\bf Case 2:} $p\in\Q_+$.
\vskip2pt

First, if $k=0$ or $N=0$, then this result can be proved as in Case~1.

Next, we assume that $k\ge 1$ and $N\ge 1$. Since $p\in\Q_+$, there exist infinitely many positive integer pairs
$(i,m)$ such that $i-pm=0$. If there are no such pairs $(i,m)$ such that $J_{i,m}\in{\cal{G}}_{k,N}$, then this result can be proved as in Case~1.
Next, we assume that there exists at least one such a pair $(i,m)$ such that $J_{i,m}\in{\cal{G}}_{k,N}$.
Assume that %$i_0$ (or equivalently, $m_0$) is the biggest such a positive integer.
$$
i_0={\rm{max}}\{i\,|\,i-pm=0 \mbox{\ and \ } J_{i,m}\in{\cal{G}}_{k,N}\}.
$$
Let $m_0=\frac{1}{p}i_0$. Then
$$
m_0={\rm{max}}\{m\,|\,i-pm=0 \mbox{\ and \ } J_{i,m}\in{\cal{G}}_{k,N}\}.
$$
We have the following three claims.

\begin{clai}
If $i_0<k$, then $I^{(1)}_{k,N}$ acts trivially on $V$, where $I^{(1)}_{k,N}=\sp_{\C}\{J_{k,m}\in{\cal{G}}_{k,N}\,|\,0\le m\le N\}$.
In particular, $V$ can be viewed as a nontrivial finite dimensional irreducible module over ${\cal{G}}_{k-1,N}$.
\end{clai}

Assume that $I^{(1)}_{k,N}$ acts nontrivially on $V$.
By the irreducibility of $V$, we have
\begin{equation}\label{claim-1}
V=I^{(1)}_{k,N}V,
\end{equation}
since $I^{(1)}_{k,N}$ is an ideal of ${\cal{G}}_{k,N}$.
%Recall an useful Lie algebra result (see, e.g., \cite{CK,CL}). %which will also be used in the following two claims.
%Let $\frak{g}$ be a finite-dimensional $\Z$-graded Lie algebra of finite depth.
%For two elements $x,y\in\frak{g}$, if $x$ is a grading operator with respect to the gradation of $\frak{g}$ and $[x,y]=y$,
%then $y$ acts nilpotently on any finite dimensional module over $\frak{g}$.
Note that $[J_{0,0},J_{k,m}]=(k-pm)J_{k,m}$, and $k-pm\ne 0$ by the assumption $i_0<k$ and the maximality of $i_0$.
Hence, each $J_{k,m}\in I^{(1)}_{k,N}$ acts nilpotently on $V$.
Since $I^{(1)}_{k,N}$ is abelian, we have $(I^{(1)}_{k,N})^{n} V=0$ for $n\gg 0$, which contradicts to \eqref{claim-1}.

\begin{clai}
If $m_0<N$, then $I^{(2)}_{k,N}$ acts trivially on $V$, where $I^{(2)}_{k,N}=\sp_{\C}\{J_{i,N}\in{\cal{G}}_{k,N}\,|\,0\le i\le k\}$.
In particular, $V$ can be viewed as a nontrivial finite dimensional irreducible module over ${\cal{G}}_{k,N-1}$.
\end{clai}

Assume that $I^{(2)}_{k,N}$ acts nontrivially on $V$. As in Claim~1, by the irreducibility of $V$, we have
\begin{equation}\label{claim-2}
V=I^{(2)}_{k,N}V,
\end{equation}
since $I^{(2)}_{k,N}$ is an ideal of ${\cal{G}}_{k,N}$.
Note that $[J_{0,0},J_{i,N}]=(i-pN)J_{i,N}$, and $i-pN\ne 0$ by the assumption $m_0<N$ and the maximality of $m_0$.
Hence, each $J_{i,N}\in I^{(2)}_{k,N}$ acts nilpotently on $V$.
Since $I^{(2)}_{k,N}$ is abelian, we have $(I^{(2)}_{k,N})^{n} V=0$ for $n\gg 0$, which contradicts to \eqref{claim-2}.

\begin{clai}
If $i_0=k$ and $m_0=N$, then $I^{(3)}_{k,N}$ acts trivially on $V$, where $I^{(3)}_{k,N}=\sp_{\C}\{J_{k,m},\,J_{i,N}\in{\cal{G}}_{k,N}\,|\,0\le m\le N,\,0\le i\le k-1\}$.
In particular, $V$ can be viewed as a nontrivial finite dimensional irreducible module over ${\cal{G}}_{k-1,N-1}$.
\end{clai}

Assume that $I^{(3)}_{k,N}$ acts nontrivially on $V$. By the irreducibility of $V$, we have
\begin{equation}\label{claim-3}
V=I^{(3)}_{k,N}V,
\end{equation}
since $I^{(3)}_{k,N}$ is an ideal of ${\cal{G}}_{k,N}$.
Consider the following decomposition of $I^{(3)}_{k,N}$:
$$
I^{(3)}_{k,N}=\C J_{i_0,m_0}+{\cal{M}}, \mbox{\ \ where\ \ }{\cal{M}}=I^{(3)}_{k,N}\backslash\C J_{i_0,m_0}.
$$
As in Claims~1 and 2, one can easily show that all elements in ${\cal{M}}$ act nilpotently on $V$.
Note that ${\cal{M}}$ is almost abelian, except $[J_{i_0,0},J_{0,m_0}]=b J_{i_0,m_0}$, where $b=-(i_0+p)m_0-i_0<0$.
To show that $I^{(3)}_{k,N}$ acts nilpotently on $V$ and then derive a contradiction to \eqref{claim-3}, we only need to show that $J_{i_0,m_0}$ acts trivially on $V$.
Note first that $J_{i_0,m_0}$ must act as a scalar $c$ on $V$, since it is a central element in ${\cal{G}}_{k,N}$.
Let $d=\dim V$ and let $Y$ be a basis of $V$. Assume that
$$
J_{i_0,0}Y=Y A,\quad J_{0,m_0} Y=Y B,\quad J_{i_0,m_0} Y=Y C,
$$
where $A$, $B$, $C$ are $d\times d$ matrices.
Then, $C=c I_d$, where $I_d$ denotes the identity matrix of order $d$.
Applying the relation $[J_{i_0,0},J_{0,m_0}]=b J_{i_0,m_0}$ to $Y$, we obtain
\begin{equation}\label{claim-3-matrix}
AB-BA=b c I_d.
\end{equation}
Considering the traces of the two sides of \eqref{claim-3-matrix}, we must have that $c=0$ (note that $b\ne0$). Hence, $J_{i_0,m_0}$ acts trivially on $V$,
and so Claim~3 holds.

Now, by Claims~1--3 and simultaneous induction on $k$ and $N$, we must have $\dim V=1$.
\QED

\vspace{12pt}

Our third lemma is the following result \cite{CK}, which will also play an important role in our classification.

\begin{lemm}\label{key-lemma-2}
Let $\cal{L}$ be a Lie superalgebra with a descending sequence of subspaces
${\cal{L}}\supset{\cal{L}}_{0}\supset{\cal{L}}_{1}\supset\ldots$ and an element $T$ satisfying $[T,{\cal L}_n] = {\cal L}_{n-1}$ for
$n\ge 1$. Let $V$ be an $\cal{L}$-module and let
$$
V_n = \{v\in V\,|\, {\cal{L}}_n v = 0\},  \ \ \ n\in\Z_+.
$$
Suppose that $V_n\ne 0$ for $n\gg 0$, and that the minimal $N\in\Z_+$ for which $V_N\ne 0$ is positive.
Then $\C[T]V_N=\C[T]\otimes_{\,\C} V_N$.
%Let $v_1, v_2, \ldots$ be a $\C$-linearly independent set of vectors of $V_N$
%generating the $\C[\pa]$-module $\C[\pa]V_N$. Then $v_1, v_2, \ldots$ is a $\C$-basis of $V_N$ and a free
%set of generators of the $\C[\pa]$-module $\C[\pa]V_N$.
In particular, $V_N$ is finite-dimensional if $V$ is a finitely generated $\C[T]$-module.
\end{lemm}

\section{The rank one case}

In this section, we shall classify all the free conformal modules of rank one over $\CB$. %without the irreducibility assumption.

Let us first construct some such conformal $\CB$-modules.
Recall that \cite{CK} a nontrivial free conformal module of rank one over $\CVir=\C[\partial]L_0$ is isomorphic to $M_{\D,\a}=\C[\partial]v$
defined by (cf.~\eqref{CVir-actions} and note that $L_0=p L$)
%with $\lambda$-actions (c.f.~\eqref{equ-CVir})
\begin{equation}\label{CVir-actions-new}
L_0\,{}_\lambda\, v=p(\partial+\D\lambda+\a)v.
\end{equation}
for some $\D,\a\in\C$.
Obviously, $M_{\D,\a}$ is also a free conformal $\CB$-module of rank one %(still denoted by $M_{\D,\a}$)
by extending the $\lambda$-actions of $L_i$ with $i\ge 1$ trivially, namely
\begin{equation}\label{trivial-extention}
L_i\,{}_\lambda\, v=
\left\{\begin{array}{ll}
p(\partial+\D\lambda+\a)v&\mbox{if \ }i=0,\\[6pt]
0&\mbox{if \ }i\ge 1.\end{array}\right.
\end{equation}
%\begin{equation}\label{trivial-extention}
%L_i\,{}_\lambda\, v=0 \ \ \ \mbox{for}\ \ \ i\ge 1.
%\end{equation}
We still denote this conformal $\CB$-module by $M_{\D,\a}$.
Clearly, $M_{\D,\a}$ is irreducible if and only if $\D\ne 0$.
The module $M_{0,\a}$ contains a unique nontrivial submodule $(\pa+ \a)M_{0,\a}$ isomorphic to $M_{1,\a}$.
We are more interested in nontrivial extensions of the conformal $\CVir$-module $M_{\D,\a}$. Let us consider $\CBB(-1)$.
For any $\b\in\C$, by replacing the $\lambda$-actions \eqref{trivial-extention} by
\begin{equation}\label{nontrivial-extention-HV}
L_i\,{}_\lambda\, v=
\left\{\begin{array}{ll}
-(\partial+\D\lambda+\a)v&\mbox{if \ }i=0,\\[6pt]
\b v&\mbox{if \ }i=1,\\[6pt]
0&\mbox{if \ }i\ge 2,\end{array}\right.
\end{equation}
we obtain a free conformal $\CBB(-1)$-module of rank one, which is a nontrivial extension of $M_{\D,\a}$ if $\b\ne 0$.
We denote this conformal $\CBB(-1)$-module by $M_{\D,\a,\b}$.
Similarly, $M_{\D,\a,\b}$ is irreducible if and only if $\D\ne 0$ or $\b\ne 0$.
The module $M_{0,\a,0}$ contains a unique nontrivial submodule $(\pa+ \a)M_{0,\a,0}$ isomorphic to $M_{1,\a,0}$.
Our main result in this section is as follows.

\begin{theo}\label{thm-rank-one-modules}
Let $M$ be a nontrivial free conformal module of rank one over $\CB$.
\baselineskip1pt\lineskip7pt\parskip-1pt
\begin{itemize}\parskip-1pt
  \item[{\rm(1)}] If $p\ne -1$, then $M\cong M_{\D,\a}$ defined by \eqref{trivial-extention} for some $\D,\a\in\C$.
  \item[{\rm(2)}] If $p=-1$, then $M\cong M_{\D,\a,\b}$ defined by \eqref{nontrivial-extention-HV} for some $\D,\a,\b\in\C$.
\end{itemize}
Furthermore, $M_{\D,\a}$ $(resp.,~M_{\D,\a,\b})$ is irreducible if and only if $\D\ne 0$ $(resp.,~\D\ne 0$ or $\b\ne0)$.
The module $M_{0,\a}$ $(resp.,~M_{0,\a,0})$ contains a unique nontrivial submodule $(\pa+ \a)M_{0,\a}$ $(resp.,~(\pa+ \a)M_{0,\a,0})$ isomorphic to $M_{1,\a}$ $(resp.,~M_{1,\a,0})$.
\end{theo}

\ni{\it Proof.}\ \
Write $M=\C[\partial]v$. First, regarding $M$ as a conformal module over $\CVir$, by \eqref{trivial-extention}, we know that
$L_0\,{}_\lambda\,v=p(\partial+\D\lambda+\a)v$ for some $\D,\a\in\C$.
By Lemma~\ref{key-lemma-1}, $L_i\,{}_\lambda\, v=0$ for all $i\gg 0$.
Assume that $k\in\Z_+$ is the largest integer such that $L_k\,{}_\lambda\, v\ne 0$.

If $k=0$, then $M$ is simply a conformal $\CVir$-module. In our notations here, $M\cong M_{\D,\a}$ if $p\ne -1$, or $M\cong M_{\D,\a,0}$ if $p=-1$.

Next, we always assume $k>0$. By the assumption $L_k\,{}_\lambda\, v\ne 0$, we can write $L_k\,{}_\lambda\, v=f(\partial,\lambda)v$, where $0\ne f(\partial,\lambda)\in\C[\partial,\lambda]$.
By relation $[L_k\,{}_\lambda\, L_k]_{\lambda+\mu}\,v=0$, we obtain
$$
f(\partial,\lambda)f(\partial+\lambda,\mu)=f(\partial,\mu)f(\partial+\mu,\lambda).
$$
Comparing the coefficients of $\lambda$, we see that $f(\partial,\lambda)$ is independent of the variable $\partial$,
and so we can denote $f(\lambda)=f(\partial,\lambda)$.
Then, by relation $[L_0\,{}_\lambda\, L_k]_{\lambda+\mu}\,v=((k+p)\lambda-p\mu)L_k\,{}_{\lambda+\mu}\,v$, we obtain
\begin{equation}\label{equ-rank-one}
(p\mu-(k+p)\lambda)f(\lambda+\mu)=p\mu f(\mu).
\end{equation}

If $k\ne -p$, then $k+p\ne 0$. By \eqref{equ-rank-one} with $\mu=0$, we immediately obtain $f(\lambda)=0$, a contradiction.

If $k=-p$ (note that this case can only occur when $p$ is a negative integer), then by \eqref{equ-rank-one} we obtain $f(\lambda+\mu)=f(\mu)$, which implies that $f(\lambda)$
is independent of the variable $\lambda$, and so we can denote $\b=f(\lambda)$.
If $p=-1$, then $k=1$, and so $M\cong M_{\D,\a,\b}$. If $p\le-2$, then $k\ge2$.
Using similar arguments as in case $k\ne -p$, one can first show that $L_i\,{}_\lambda\,v=0$ for $1\le i\le k-1$.
Then, by relation $[L_1\,{}_\lambda\, L_{k-1}]_{\lambda+\mu}\,v=0$, we obtain $(\lambda+(1+p)\mu)\b=0$, which implies $\b=0$, a contradiction.
This completes the proof.
\QED

\section{Classification theorem}

Now, we can show that the irreducible modules appeared in Theorem~\ref{thm-rank-one-modules} exhaust
all nontrivial finite irreducible conformal modules over $\CB$. Namely, we have the following classification.

\begin{theo}\label{thm-classification}
Let $M$ be a nontrivial finite irreducible conformal module over $\CB$.
\baselineskip1pt\lineskip7pt\parskip-1pt
\begin{itemize}\parskip-1pt
  \item[{\rm(1)}] If $p\ne -1$, then $M\cong M_{\D,\a}$ defined by \eqref{trivial-extention} for some $\D,\a\in\C$ with $\D\ne 0$.
  \item[{\rm(2)}] If $p=-1$, then $M\cong M_{\D,\a,\b}$ defined by \eqref{nontrivial-extention-HV} for some $\D,\a,\b\in\C$ with $\D\ne 0$ or $\b\ne0$.
\end{itemize}
\end{theo}

Let $M$ be a nontrivial finite irreducible conformal module over $\CB$.
Our basic strategy is to show that $M$ must be free of rank one (see Lemma~\ref{free-of-rank-one} below) by using the three technical lemmas given in Section~3.
Then, the above result will follow from Theorem~\ref{thm-rank-one-modules}.

\begin{lemm}\label{free-of-rank-one}
The conformal $\CB$-module $M$ must be free of rank one.
\end{lemm}
\ni{\it Proof.}\ \
First, by Lemma~\ref{key-lemma-1}, the $\lambda$-action of $L_i\in\CB$ on $M$ is trivial for all $i\gg 0$.
Assume that $k\in\Z_+$ is the largest integer such that the $\lambda$-action of $L_k$ on $M$ is nontrivial.
Then $M$ is simply a nontrivial finite irreducible conformal module over $\CB_{[k]}$, where $\CB_{[k]}$ is defined by \eqref{equ-quotiont}.
%Recall that $\CB_{[0]}\cong\CVir$, by \eqref{CVir-actions}, we may assume that $k\ge 1$.

Furthermore, by Proposition~\ref{prop-observation}, the conformal $\CB_{[k]}$-module $M$
can be viewed as a module over the associated extended annihilation algebra
${\cal A}(\CB_{[k]})^e$ satisfying
\begin{equation}\label{equiverlent}
\bar{L}_{i,m} v = 0\ \ \ \mbox{for}\ \ \ v\in M,\ \ \ 0\le i\le k,\ \ \ m\gg 0.
\end{equation}
%Note that $T-\frac{1}{p} \bar{L}_{0,-1}$ is a central element of ${\cal A}(\CB_{[k]})^e$ by Remark~\ref{remark-quotionts},
%we see that ${\cal A}(\CB_{[k]})^e$ is a direct sum as ideals of the commutative
%Lie algebra $\C(T-\frac{1}{p} \bar{L}_{0,-1})$ and the Lie algebra ${\cal A}(\CB_{[k]})$.
%%Namely,
%%$$
%%{\cal A}(\CB_{[k]})^e=\C(T-\frac{1}{p} \bar{L}_{0,-1})\oplus {\cal A}(\CB_{[k]}).
%%$$
%Hence, $T-\frac{1}{p} \bar{L}_{0,-1}$ acts as a scalar, and ${\cal A}(\CB_{[k]})$ acts irreducibly on $M$.

For simplicity, in what follows, we denote ${\cal{L}}={\cal A}(\CB_{[k]})^e$.
Set
$$
{\cal{L}}_n=\sp_{\C}\{\bar{L}_{i,m}\in {\cal{L}}\,|\,0\le i\le k,\,m\ge n-1\}, \ \ \  n\in\Z_+.
$$
Note that ${\cal{L}}_{0}={\cal A}(\CB_{[k]})$.
Clearly, ${\cal{L}}\supset{\cal{L}}_{0}\supset{\cal{L}}_{1}\supset\ldots$ and, by Remark~\ref{remark-quotionts},
the element ${T}\in{\cal{L}}$ satisfies $[{T},{\cal L}_n] = {\cal L}_{n-1}$ for $n\ge 1$.
For the ${\cal{L}}$-module $M$, we set
%$$
%M_n = \{v\in M\,|\, {\cal{L}}_n v = 0\}=\{v\in M\,|\, \bar{L}_{i,m} v = 0\ \mbox{for}\ 0\le i\le k,\,m\ge n-1\},  \ \ \ n\in\Z_+.
%$$
$$
M_n = \{v\in M\,|\, {\cal{L}}_n v = 0\},  \ \ \ n\in\Z_+.
$$
By \eqref{equiverlent}, $M_n\ne 0$ for $n\gg 0$. Assume that $N\in\Z_+$ is the smallest integer such that $M_N\ne 0$.

Suppose $N=0$. Take $0\ne v\in M_0$. Then $U({\cal{L}})v=\C[{T}]U({\cal{L}}_0)v=\C[{T}]v$,
and so $M=\C[{T}]v$ by the irreducibility of $M$. Since ${\cal{L}}_0$ is an ideal of ${\cal{L}}$,
we see that ${\cal{L}}_0$ acts trivially on $M$. Hence, $M$ is simply an irreducible $\C[{T}]$-module,
and so $M$ is one-dimensional. Equivalently, $M$ is a one-dimensional trivial conformal $\CB$-module, a contradiction.
%If $N=0$, then it is easy to see that $M$ is a one-dimensional trivial ${\cal{L}}$-module, a contradiction.

Next, we always assume $N\ge 1$.
By Remark~\ref{remark-quotionts}, $T-\frac{1}{p} \bar{L}_{0,-1}$ is a central element in $\cal{L}$.
Hence, $T-\frac{1}{p} \bar{L}_{0,-1}$ acts on $M$ as a scalar, and ${\cal{L}}_{0}$ acts irreducibly on $M$.
Furthermore,  by relation $\bar{L}_{i,-1}=\frac{1}{p}[\bar{L}_{i,0},\bar{L}_{0,-1}]$, we see that
the action of ${\cal{L}}_{0}$ is determined by ${\cal{L}}_{1}$ and $\bar{L}_{0,-1}$ (or equivalently, determined by ${\cal{L}}_{1}$ and $T$).
Clearly, $M_N$ is ${\cal{L}}_1$-invariant.
By the irreducibility of $M$ and Lemma~\ref{key-lemma-2}, we see that $M=\C[T]\otimes_{\,\C} M_N$ and  $M_N$ is a nontrivial irreducible finite-dimensional ${\cal{L}}_1$-module.

If $N=1$, then by the definition of $M_1$, we see that $M_1$ is a trivial ${\cal{L}}_1$-module, a contradiction.

If $N\ge 2$, then by the definition of $M_N$, we see that $M_N$ can be viewed as a ${\cal{L}}_1/{\cal{L}}_N$-module.
Note that  ${\cal{L}}_1/{\cal{L}}_N\cong {\cal G}_{k,N-2}$.  By Lemma~\ref{key-lemma-3}, we must have that
$M_N$ is one-dimensional. Equivalently, $M$ is free of rank one as a conformal $\CB$-module.
\QED

\vspace{12pt}

At last, by Theorems~\ref{thm-rank-one-modules} and \ref{thm-classification},
one can easily obtain the classification of finite irreducible conformal modules over the finite Lie conformal algebra ${\frak b}(n)$ for $n\ge 1$.
Recall that ${\frak b}(n)$ has a $\C[\partial]$-basis $\{\bar{L}_i\,|\,0\le i\le n\}$
with $\lambda$-brackets (cf.~\eqref{brackets} and \eqref{equ-quotiont-special})
\begin{equation*}
[\bar{L}_i\, {}_\lambda \, \bar{L}_j] =
\left\{\begin{array}{ll}
((i-n)\partial+(i+j-2n)\lambda) \bar{L}_{i+j}&\mbox{if \ }i+j\le n,\\[6pt]
0&\mbox{otherwise,}\end{array}\right.
\end{equation*}
and  contains a conformal subalgebra $\C[\partial]\bar{L}_0\cong\CVir$, which has a free conformal module of rank one $M_{\D,\a}$ given by \eqref{CVir-actions-new} with $p=-n$.
Clearly, there is a free conformal ${\frak b}(n)$-module of rank one $\C[\partial]v$ defined by
\begin{equation}\label{bn-trivial-extention}
\bar{L}_i\,{}_\lambda\, v=
\left\{\begin{array}{ll}
-n(\partial+\D\lambda+\a)v&\mbox{if \ }i=0,\\[6pt]
0&\mbox{if \ }1\le i\le n,\end{array}\right.
\end{equation}
where $\D,\a\in\C$. We still denote this conformal ${\frak b}(n)$-module by $M_{\D,\a}$.
In addition, consider the special case ${\frak b}(1)$.  Replacing the above $\lambda$-actions by
\begin{equation}\label{bn-nontrivial-extention-HV}
\bar{L}_i\,{}_\lambda\, v=
\left\{\begin{array}{ll}
-(\partial+\D\lambda+\a)v&\mbox{if \ }i=0,\\[6pt]
\b v&\mbox{if \ }i=1,\end{array}\right.
\end{equation}
where $\b\in\C$,
we obtain a free conformal ${\frak b}(1)$-module of rank one, which is a nontrivial extension of the conformal $\CVir$-module $M_{\D,\a}$ if $\b\ne 0$.
We still denote this conformal ${\frak b}(1)$-module by $M_{\D,\a,\b}$.
Since ${\frak b}(n)$ is a quotient algebra of $\CBB(-n)$, by Theorems~\ref{thm-rank-one-modules} and \ref{thm-classification}, we have

\begin{coro}\label{classification-b}
Let $M$ be a nontrivial free conformal module of rank one over ${\frak b}(n)$.
\baselineskip1pt\lineskip7pt\parskip-1pt
\begin{itemize}\parskip-1pt
  \item[{\rm(1)}] If $n>1$, then $M\cong M_{\D,\a}$ defined by \eqref{bn-trivial-extention} for some $\D,\a\in\C$.
  \item[{\rm(2)}] If $n=1$, then $M\cong M_{\D,\a,\b}$ defined by \eqref{bn-nontrivial-extention-HV} for some $\D,\a,\b\in\C$.
\end{itemize}
Furthermore, $M_{\D,\a}$ $(resp.,~M_{\D,\a,\b})$ is irreducible if and only if $\D\ne 0$ $(resp.,~\D\ne 0$ or $\b\ne0)$.
The module $M_{0,\a}$ $(resp.,~M_{0,\a,0})$ contains a unique nontrivial submodule $(\pa+ \a)M_{0,\a}$ $(resp.,~(\pa+ \a)M_{0,\a,0})$ isomorphic to $M_{1,\a}$ $(resp.,~M_{1,\a,0})$.
The modules $M_{\D,\a}$ with $\D\ne 0$ and  $M_{\D,\a,\b}$ with $\D\ne 0$ or $\b\ne0$
exhaust all nontrivial finite irreducible conformal modules over ${\frak b}(n)$.
\end{coro}

\vskip10pt

\small \ni{\bf Acknowledgement}
This work was supported by the National Natural Science Foundation of China (Nos.~11371278, 11431010, 11401570, 11661063, 11301109),
and the Natural Science Foundation of Jiangsu Province, China (No.~BK20140177).
%and the Fundamental Research Funds for the Central Universities (No.~2014QNA68).


\begin{thebibliography}{9999}\vskip0pt\small\footnotesize
\def\re{\bibitem}\parindent=2ex\parskip=-8pt\baselineskip=-2pt

\bibitem{ADKP} E. Arbarello, C. De Concini, V. Kac, C. Procesi,
Moduli spaces of curves and representation theory,
{\it Comm. Math. Phys.} {\bf 117} (1998) 1--36.

\bibitem{BKV} B. Bakalov, V. Kac, A. Voronov,
Cohomology of conformal algebras,
{\it Comm. Math. Phys.} {\bf 200} (1999) 561--598.

\bibitem{BKL1}  C. Boyallian, V. Kac, J. Liberati,
On the classification of subalgebras of $Cend_N$ and $gc_N$,
{\it J. Algebra} {\bf 260} (2003) 32--63.

\bibitem{BKL2}  C. Boyallian, V. Kac, J. Liberati,
Finite growth representations of infinite Lie conformal algebras,
{\it J. Math. Phys.} {\bf 44} (2) (2003) 754--770.

\bibitem{BKL3}  C. Boyallian, V. Kac, J. Liberati,
Classification of finite irreducible modules over the Lie conformal superalgebra $CK_6$,
{\it Comm. Math. Phys.} {\bf 317} (2) (2013) 503--546.

\bibitem{BM} C. Boyallian, V. Meinardi,
Finite growth representations of conformal Lie algebras that contain a Virasoro subalgebra,
{\it J. Algebra} {\bf 388} (2013) 141--159.

\bibitem{CK} S. Cheng, V. Kac, Conformal modules,
{\it Asian J. Math.} {\bf 1} (1) (1997) 181--193, {\it Asian J. Math.}  {\bf 2} (1) (1998) 153--156 (Erratum).

\bibitem{CKW}  S. Cheng, V. Kac, M. Wakimoto,
Extensions of conformal modules, in: Topological Field Theory, Primitive Forms and Related
Topics, Proceedings of Taniguchi Symposium, in: Progr. Math., vol. 160, Birkh\"{a}user, 1998.

%\bibitem{CKW} S. Cheng, V. Kac, M.Wakimoto, ``Extensions of Neveu-Schwarz modules,'' {\it J. Math. Phys}. {\bf 41}, 2271--2294 (2000).

%\bibitem{CL} S. Cheng, N. Lam, Finite conformal modules over $N=2,3,4$ superconformal algebras, {\it J. Math. Phys}. {\bf 42}, 906--933 (2001).

\bibitem{DK} A. D'Andrea, V. Kac,
Structure theory of finite conformal algebras,
{\it Sel. Math.} (N.S.) {\bf 4} (3) (1998) 377--418.

\bibitem{DeK1} A. De Sole, V. Kac,
Subalgebras of $gc_N$ and Jacobi polynomials,
{\it Canad. Math. Bull.} {\bf 45} (4) (2002) 567--605.

\bibitem{DeK2} A. De Sole, V. Kac,
Lie conformal algebra cohomology and the variational complex,
{\it Comm. Math. Phys.} {\bf 292} (2009) 667--719.

\bibitem{FSW} G. Fan, Y. Su, H. Wu,
Loop Heisenberg-Virasoro Lie conformal algebra,
{\it J. Math. Phys.} {\bf 55} (2014) 123508.

\bibitem{FSX} G. Fan, Y. Su, C. Xia,
Infinite rank Schr$\ddot{o}$dinger-Virasoro type Lie conformal algebras,
{\it J. Math. Phys.} {\bf 57} (2016) 081701.

\bibitem{FK} D. Fattori, V. Kac,
Classification of finite simple Lie conformal superalgebras,
in: Special Issue in Celebration of Claudio Procesi's 60th Birthday,
{\it J. Algebra} {\bf 258} (1) (2002) 23--59.

\bibitem{FKR} D. Fattori, V. Kac, A. Retakh,
Structure theory of finite Lie conformal superalgebras, in: Lie Theory and Its Applications in
Physics V, World Sci. Publ., River Edge, NJ, 2004, pp. 27--63.

\bibitem{He} M. Henkel,
Schr$\ddot{o}$dinger invariance and strongly anisotropic critical systems,
{\it J. Stat. Phys.} {\bf 75} (5) (1994) 1023--1029.

\bibitem{K1} V. Kac, {\it Vertex Algebras for Beginners},
Univ. Lecture Ser., vol. 10, Amer. Math. Soc., 1998.

\bibitem{K2} V. Kac,
Formal distribution algebras and conformal algebras,
A talk at the Brisbane Congress in Math. Physics, 1997.

\bibitem{K3} V. Kac, The idea of locality, in: H.-D. Doebner, et al. (Eds.),
Physical Applications and Mathematical Aspects of Geometry,
Groups and Algebras, World Sci. Publ., Singapore, 1997, pp. 16--32.

\bibitem{MZ}  C. Mart\'{\i}nez, E. Zelmanov,
Irreducible representations of the exceptional Cheng-Kac superalgebra,
{\it Trans. Amer. Math. Soc.} {\bf 366} (11) (2014) 5853--5876.

%\bibitem{RU} C. Roger, J. Unterberger, The Schr$\ddot{o}$dinger-Virasoro Lie group and algebra: representation theory and cohomological study,
%{\it Ann. Henri Poincar$\acute{e}$} {\bf 7}, (2006) 1477--1529.

\bibitem{S1} Y. Su,
Low dimensional cohomology of general conformal algebras $gc_N$,
{\it J. Math. Phys.} {\bf 45} (2004) 509--524.

%\bibitem{S2} Y. Su, Quasifinite representations of a Lie algebra of Block  type,
%\textit{J. Algebra} \textbf{276} (2004) 117--128.

\bibitem{SXX1} Y. Su, C. Xia, Y. Xu,
Quasifinite representations of a class of Block type Lie algebras $\BB$,
\textit{J. Pure Appl. Algebra} \textbf{216} (4) (2012) 923--934.

\bibitem{SXX2} Y. Su, C. Xia, Y. Xu,
Classification of quasifinite representations of a Lie algebra related to Block type,
\textit{J. Algebra} \textbf{393} (2013) 71--78.

\bibitem{SY} Y. Su, L. Yuan,
Schr$\ddot{\rm o}$dinger-Virasoro Lie conformal algebra,
{\it J. Math. Phys.} {\bf 54} (2013) 053503.

\bibitem{SYue} Y. Su, X. Yue,
Filtered Lie conformal algebras whose associated graded algebras are isomorphic to that of general conformal algebra $gc_1$,
{\it J. Algebra} {\bf 340} (2011) 182--198.

\bibitem{WCY} H. Wu, Q. Chen, X. Yue,
Loop Virasoro Lie conformal algbera,
{\it J. Math. Phys.} {\bf 55} (2014) 011706.

\bibitem{XZ} C. Xia, R. Zhang,
Unitary highest weight modules over Block type Lie algebras $\BB$,
\textit{J. Lie theory} {\bf 23} (1) (2013) 159--176.

\bibitem{Z1} E. Zelmanov,
On the structure of conformal algebras, in: Combinatorial and Computational Algebra, Hong Kong, 1999, in:
Contemp. Math., vol. 264, Amer. Math. Soc., Providence, RI, 2000, pp. 139--153.

\bibitem{Z2} E. Zelmanov,
Idempotents in conformal algebras, in: Proceedings of the Third International Algebra Conference, Tainan,
2002, Kluwer Acad. Publ., Dordrecht, 2003, pp. 257--266.
\end{thebibliography}
\end{document}